\documentclass[a4paper, 12pt]{extarticle}
\usepackage{mathtext}
\usepackage[T2A]{fontenc}
\usepackage[english]{babel}
\usepackage{latexsym,amsfonts,amssymb,amsmath,longtable,amsthm,setspace}
\usepackage[pic]{xy}
\usepackage{bbm,makeidx}
\usepackage[centerlast,small]{caption}
\usepackage{amsmath}
\usepackage[cp1251]{inputenc}
\def\G{{\Gamma}}

\def\Z{{\mathbb Z}}

\def\cS{{\mathfrak S}}
\def\Aut{\operatorname{Aut}}
\def\Out{\operatorname{Out}}
\def\Id{\operatorname{Id}}
\def\Inn{\operatorname{Inn}}
\def\Ind{\operatorname{Ind}}
\newcommand{\Mat}[4]{\left( \begin{array}{cc}
                            #1 & #2 \\
                            #3 & #4
                      \end{array} \right)}

\renewcommand{\ge}{\geqslant}

\newtheorem{lem}{{\scshape Lemma}}

\newtheorem{ttt}{{\scshape Theorem}}
\newtheorem{prp}{{\scshape Proposition}}

\begin{document}
\title{The $R_{\infty}$  and $S_{\infty}$ properties for linear algebraic groups}
\author{Alexander Fel'shtyn}

\author{Alexander Fel'shtyn and Timur Nasybullov\footnote{The author is supported by  Russian Science Foundation (project 14-21-00065)}}

\maketitle

\begin{abstract}
In this paper we study twisted conjugacy classes and isogredience classes for automorphisms of reductive linear algebraic groups. We show that reductive linear algebraic groups over some fields of zero characteristic possess the  $R_\infty$ and $S_\infty$ properties.

~
\end{abstract}

\section{Introduction}

Let $\varphi:G\to G$ be an endomorphism
of a  group $G$.
Then two elements
$x,y$ of $G$ are said to be twisted $\varphi$-conjugate, if there exists a third element
$z\in G$ such that $x = z y \varphi(z)^{-1}$. The equivalence classes  are called the twisted conjugacy classes or the  Reidemeister classes  of $\varphi$. The Reidemeister number of $\varphi$ denoted by $R(\varphi)$, is the number of those twisted conjugacy classes of $\varphi$. This number is either a positive integer or $\infty$ and we do not distinguish different infinite cardinal numbers.
An infinite group $G$ has the $R_\infty$-property if for every \underline{automorphism} $\varphi$ of $G$ the Reidemeister number of $\varphi$ is
infinite.

The interest in twisted conjugacy relations has its origins, in
particular, in the Nielsen-Reidemeister fixed point theory (see,
e.g. \cite{FelshB, Jiang}), in Arthur- Selberg theory (see, eg.
\cite{Arthur, Shokra}), in  algebraic geometry (see, e.g.
\cite{Groth}), in Galois cohomology \cite{Serre} and in the theory of linear algebraic groups (see, e.g. \cite{S}). In representation theory twisted conjugacy probably
occurs first in Gantmacher's paper \cite{Gantmacher} (see, e.g
\cite{Springer,Onishik-Vinberg})

The problem of determining which classes of discrete infinite groups have the-$R_{\infty}$ property is an area of active research initiated by Fel'shtyn and Hill in 1994 \cite{FelHill}.
Later, it was shown by various authors that the following groups
have the $R_\infty$-property:
 non-elementary Gromov hyperbolic groups
\cite{FelPOMI,ll}; relatively hyperbolic groups \cite{f07};
 Baumslag-Solitar groups $BS(m,n)$
except for $BS(1,1)$ \cite{FelGon08Progress}, generalized
Baumslag-Solitar groups, that is, finitely generated groups which
act on a tree with all edge and vertex stabilizers infinite cyclic
\cite{LevittBaums}; the solvable generalization $\G$ of $BS(1,n)$
given by the short exact sequence $1 \rightarrow \mathbb
Z[\frac{1}{n}] \rightarrow \G \rightarrow \mathbb Z^k \rightarrow
1,$ as well as any group quasi-isometric to $\G$ \cite{TabWong};
a wide class of saturated weakly branch groups (including the
Grigorchuk group \cite{GrFA} and the Gupta-Sidki group
\cite{GuSi}) \cite{FelLeonTro}, Thompson's groups $F$ \cite{bfg} and $T$ \cite{BMV, GoSan};
generalized Thompson's groups $F_{n,\:0}$ and their finite direct
products \cite{GK}; Houghton's groups \cite{DS, JLL}; symplectic groups ${\rm Sp}(2n,\mathbb Z)$,
the mapping class groups $Mod_{S}$ of a compact oriented surface
$S$ with genus $g$ and $p$ boundary components, $3g+p-4>0$, and
the full braid groups $B_n(S)$ on $n>3$ strands of a compact
surface $S$ in the cases where $S$ is either the compact disk $D$,
or the sphere $S^2$ \cite{dfg}; some classes of Artin groups of
infinite type \cite{juhasz};
  extensions of  ${\rm SL}(n,\Z)$, ${\rm PSL}(n,\Z)$, ${\rm GL}(n,\Z)$, ${\rm PGL}(n,\Z)$, ${\rm Sp}(2n,\Z)$,
${\rm PSp}(2n,\Z)$, $n>1$, by a countable abelian group, and normal
subgroups of ${\rm SL}(n,\Z)$, $n>2$, not contained in the center
 \cite{MS}; ${\rm GL}(n,K)$
and ${\rm SL}(n,K)$ if $n>2$ and
$K$ is an infinite integral domain
with trivial group of automorphisms, or
$K$ is an integral domain, which has
a zero characteristic and for which $\Aut(K)$ is periodic
\cite{Nasybull2012};
irreducible lattices in a connected semisimple Lie group $G$ with finite center and real rank
at least 2 \cite{MubeenaSankaran2014TrGr};
non-amenable, finitely generated residually finite groups \cite{FT}
(this class gives a lot of new examples of groups with the $R_\infty$-property); some metabelian
groups of the form ${\mathbb Q}^n\rtimes \mathbb Z$
and ${\mathbb Z[1/p]}^n\rtimes \mathbb Z$ \cite{FelGon2011Q};
lamplighter groups
$\mathbb Z_n \wr \mathbb Z$ if and only if $2|n$ or
$3|n$ \cite{gowon1};
 free nilpotent groups  $N_{rc}$ of rank $r=2$ and class $c\ge 9$
 \cite{GoWon09Crelle}, $N_{rc}$ of rank
$r = 2$ or $r = 3$ and class $c \geq 4r,$ or rank $r \geq 4$ and
class $c \geq 2r,$ any group $N_{2c}$ for $c \geq 4$, every free
solvable group $S_{2t}$ of rank 2 and class $t \geq 2$ (in
particular the free metabelian group $M_2 = S_{22}$ of rank 2),
any free solvable group $S_{rt}$ of rank $r \geq 2$ and class $t$
big enough \cite{Romankov}; some crystallographic groups
\cite{DePe2011,LutScze2011}. Recently, in \cite{DeGon} it was proven that
$N_{rc}$, $r > 1$ has the $R_\infty$-property if and only if $c \geq 2r$.

Let $\Psi$ belongs to $\Out (G)= \Aut (G)/ \Inn (G)$. We consider an outer automorphism $\Psi \in \Out (G)$ as a collection of ordinary automorphisms
$a \in  \Aut (G) $. We say that two automorphisms
$a,b \in \Psi$ are \emph{similar} (or \emph{isogredient})
if $b=\varphi_h a  \varphi_h^{-1}$ for some $h\in G$,
where $\varphi_h(g)=hgh^{-1}$ an inner automorphism induced by the element $h$ (see \cite{ll}).
Let $\cS(\Psi)$ be the set of isogredience classes of automorphisms representing $\Psi$.
Denote by $S(\Psi)$ the cardinality of the set $\cS(\Psi)$.
A group $G$ is called an $S_\infty$-\emph{group} if for
any  $\Psi$ the set $\cS(\Psi)$ is infinite,
i.~e. $S(\Psi)=\infty$ (see \cite{FTMPI}).

In this paper we study the $R_\infty$ and  $S_\infty$  properties  for   linear algebraic groups.
First results in this direction were obtained  for some classes of Chevalley groups by Nasybullov  in \cite{Nas2}.

In Section 3 we extend the previous result from \cite{Nas2} and prove

\noindent \textbf{Theorem 2.} \emph{ Let $G$ be a Chevalley group of the type $\Phi$ over the field $F$ of zero characteristic. If the transcendence degree of $F$ over $\mathbb{Q}$ is finite, then $G$ possesses the $R_{\infty}$-property.}

 The following main theorem is proved in Section 4.

 \noindent \textbf{Theorem 3.} \emph{Let $F$ be such an algebraically closed field of zero characteristic that the transcendence degree of $F$ over $\mathbb{Q}$ is finite. If the reductive linear algebraic group $G$ over the field $F$ has a nontrivial quotient group $G/R(G)$, where $R(G)$ is the radical of $G$, then $G$ possesses the $R_{\infty}$-property.
}

These theorems can not be generalized to groups over a field of non-zero characteristic. It follows from the following theorem of Steinberg \cite[Theorem 10.1]{S}.

 \noindent \textbf{Theorem.} \emph{Let $G$ be a connected linear algebraic group and $\varphi$ be an endomorphism of $G$ onto $G$. If $\varphi$ has a finite set of fixed points, then $G=\{x\varphi(x^{-1})~|~x\in G\}$.
}

We would like to point out that R.~Steinberg \cite{S} calls by an automorphism of a linear algebraic group a bijective endomorphism which is a morphism and its inverse is a morphism too. However, throughout the paper we understand an automorphism of a linear algebraic group as an automorphism of an abstract group, i.~e. a bijective endomorphism of a group.

Any semisimple linear algebraic group over an algebraically closed field of positive characteristic possesses an automorphism $\varphi$ with finitely many fixed points (Frobenius morphism, see \cite[\S 3.2]{Z}), therefore, this group coincides with the set $\{x\varphi(x^{-1})~|~x\in G\}=[e]_{\varphi}$, hence $R(\varphi)=1$ and such a group can not have the $R_{\infty}$-property.

If $T_1, T_2,\dots$ are algebraically independent over $\mathbb{Q}$ elements, then the fields $\overline{\mathbb{Q}}$, $\overline{\mathbb{Q}(T_1,\dots,T_k)}~(k\geq 1)$ are algebraically closed fields of zero characteristic with finite transcendence degree over $\mathbb{Q}$. Then the reductive linear algebraic groups over these fields possess the $R_{\infty}$-property.

In the Section 5 we prove that an infinite  reductive linear algebraic  group $G$ over the field $F$ of zero characteristic and  finite transcendence degree over $\mathbb{Q}$ which possesses an automorphism $\varphi$ with a finite Reidemeister number is a torus.

In the Section 6 we prove the following

\noindent \textbf{Theorem 5.} \emph{Let $F$ be such an algebraically closed field of zero characteristic that the transcendence degree of $F$ over $\mathbb{Q}$ is finite. If the reductive linear algebraic group $G$ over the field $F$ has a nontrivial quotient group $G/R(G)$, then $G$ possesses the $S_{\infty}$-property.
}

\bigskip
\textbf{Acknowledgment.}
The authors are grateful to Andrzej D\c abrowski,  Evgenij Troitsky and Evgeny Vdovin for the numerous important discussions on linear algebraic groups.
The first author would like to thank the Max Planck Institute for Mathematics(Bonn) for its kind support and hospitality while
a part of this work was completed.

\section{Preliminaries}
In this paragraph we recall some preliminary statements which are used in the paper.  A lot of used results are thoroughly written in \cite{Nas2}, the reader can use it as a background material.

Symbols $I_n$ and $O_{n\times m}$ mean the identity $n\times n$ matrix and the $n\times m$ matrix with zero entries, respectively. If  $A$ an $n\times n$ matrix and $B$ an $m\times m$ matrix, then the symbol  $A\oplus B$ denotes the direct sum of the matrices $A$ and $B$, i.~e. the block-diagonal $(m+n)\times (m+n)$ matrix
$$
\newcommand{\tempa}{\multicolumn{1}{c|}{A}}
\newcommand{\tempb}{\multicolumn{1}{|c}{B}}
\begin{pmatrix}
\tempa&O_{n\times m}\\\cline{1-2}
O_{m\times n}&\tempb
\end{pmatrix}
.
$$
  It is obvious that for a pair of $n\times n$ matrices $A_1, A_2$ and for a pair of $m\times m$ matrices $B_1,B_2$ we have
$(A_1\oplus B_1)(A_2\oplus B_2)=A_1A_2\oplus B_1B_2$, $(A_1\oplus B_1)^{-1}=A_1^{-1}\oplus B_1^{-1}$.

Symbols $G\times H$ and $G\circ H$ mean the direct product and the central product of the  groups $G$ and $H$, respectively.

If $g$ is an element of the group $G$, then $\varphi_g$ denotes an inner automorphism induced by the element $g$. The following lemma can be found in \cite[Corollary 2.5]{FTC}.
\begin{lem}\label{pr1} Let $\varphi$ be an automorphism of the group $G$ and $\varphi_g$ be an inner au\-to\-mor\-phism of the group $G$. Then $R(\varphi\varphi_g)=R(\varphi)$.
\end{lem}
The next lemma is proved in \cite[Lemma 2.1]{MS}
\begin{lem}\label{pr3} Let
$$1\to N\to G \to A\to 1$$
be an exact sequence of groups. Suppose that $N$ is a characteristic subgroup of $G$ and that $A$ possesses the $R_{\infty}$-property, then $G$ also possesses the $R_{\infty}$-property.
\end{lem}
Here we prove similar result for the $S_{\infty}$-property.
\begin{lem}\label{an3} Let
$$1\to N\to G \to A\to 1$$
be an exact sequence of groups. Suppose that $N$ is a characteristic subgroup of $G$ and that $A$ possesses the $S_{\infty}$-property, then $G$ also possesses the $S_{\infty}$-property.
\end{lem}
\textbf{Proof.} Let $\varphi$ be an automorphism of the group $G$. Since $N$ is a characteristic subgroup of $G$ then $\varphi$ induces an automorphism $\overline{\varphi}$ of the group $A$. Since the group $A$ has the $S_{\infty}$-property then there exists an infinite set of elements $\overline{g}_1, \overline{g}_2, \dots$ of the group $A$ such that $\varphi_{\overline{g}_i}\overline{\varphi}$ and $\varphi_{\overline{g}_j}\overline{\varphi}$ are not isogredient for $i\neq j$.

Suppose that $S(\varphi {\rm Inn}(G))<\infty$. Then there exists a pair of isogredient automorphisms in the set $\varphi_{g_1}\varphi, \varphi_{g_2}\varphi,\dots$. Suppose that $\varphi_{g_i}\varphi$ and $\varphi_{g_j}\varphi$ are isogredient for $i\neq j$. Then for some element $h\in G$ we have
$$\varphi_{g_i}\varphi=\varphi_h\varphi_{g_j}\varphi\varphi_h^{-1}.$$
From this equality we have the following equality in the group ${\rm Aut}(A)$
$$\varphi_{\overline{g}_i}\overline{\varphi}=\varphi_{\overline{h}}\varphi_{\overline{g}_j}\overline{\varphi}\varphi_{\overline{h}}^{-1},$$
but it contradicts to the choice of the elements $\overline{g}_1,\overline{g}_2,\dots$\hfill$\square$

Let $\nu$ be a map from the set of rational numbers $\mathbb{Q}$ to the set $2^{\pi}$ of all subsets of the set of prime numbers $\pi$, which acts on the irreducible fraction $x=a/b$ by the rule $$\nu(x)=\{all ~the ~prime ~devisors ~of ~a\}\cup \{all ~the ~prime ~devisors ~of ~b\}.$$

 The proof of the following lemma is presented in \cite[Lemma 5]{Nas2}.
\begin{lem}\label{lem1}\notag Let $F$ be a field of zero characteristic and $x_1,\dots,x_k$ be elements of $F$ which are algebraically independent over the field $\mathbb{Q}$. Let $x_{k+1}$ be such an element of $F$, that the elements $x_1,\dots, x_{k+1}$ are algebraically dependent over $\mathbb{Q}$. Let $\delta$ be an automorphism of the field $F$ which acts on this elements by the rule
$$
\delta: x_i\mapsto t_0 t_i x_i,~~~ i=1,\dots,k+1,
$$
where $t_0,\dots,t_{k+1}\in \mathbb{Q}$ and $t_1,\dots,t_{k+1}$ are not equal to $1$. If $\nu(t_i)\cap\nu(t_j)=\varnothing$ for $i\neq j$, then $x_{k+1}=0$.
\end{lem}

Using this lemma we prove the following auxiliary statement.
\begin{lem}\label{lem2} Let $F$ be such a field of zero characteristic, that the transcendence degree of $F$ over $\mathbb{Q}$ is finite. If the automorphism $\delta$ of the field $F$ acts on the elements $z_1$, $z_2,\dots$ of the field $F$ by the rule
$$\delta: z_i \mapsto \alpha a_iz_i,$$
where $\alpha\in F$, $1\neq a_i\in \mathbb{Q}\subseteq F$ and $\nu(a_i)\cap\nu(a_j)=\varnothing$ for $i\neq j$, then there are only a finite number of non-zero elements in the set $z_1,z_2,\dots$.
\end{lem}
\textbf{Proof.} If all the elements $z_1,z_2,\dots$ are equal to zero then there is nothing to prove. Hence we can consider that there exists a non-zero element in the set $z_1,z_2,\dots$ Without loosing of generality we can consider that $z_1\neq 0$ (Otherwise we can reenumerate the elements $z_1,z_2,\dots$ and do the first element not to be equal to zero. If the statement holds for the reenumerated set, then it holds for the original set $z_1,z_2,\dots$). Let us denote $y_i=z_iz_1^{-1}$. Then the automorphism $\delta$ acts on the element $y_i$ by the rule
$$\delta(y_i)=\delta(z_iz_1^{-1})=\delta(z_i)\delta(z_1^{-1})=\alpha a_iz_i\alpha^{-1} a_1^{-1}z_1^{-1}=a_ia_1^{-1}z_iz_1^{-1}=a_ia_1^{-1}y_i$$

Since the transcendence degree of $F$ over $\mathbb{Q}$ is finite, then there exists a maximal subset of algebraically independent over $\mathbb{Q}$ elements in the set $y_2,y_3,\dots$, i.e. there exists such a finite set  $y_{i_1},y_{i_2},\dots,y_{i_k}$ of algebraically independent over $\mathbb{Q}$ elements, that the set $y_{i_1},y_{i_2},\dots,y_{i_k},y_j$ is algebraically dependent over $\mathbb{Q}$ for every $j$.

Without loosing of generality we can consider that the set $y_2,\dots,y_k$ is a maximal subset of algebraically independent over $\mathbb{Q}$ elements in the set $y_1, y_2, \dots$

If $n>k$ is a positive integer, then the elements $y_2,\dots,y_k, y_n\in F$ satisfy the conditions of the lemma \ref{lem1}. Therefore $y_n=0$ for all $n>k$ and since $y_n=z_nz_1^{-1}$ then $z_n=0$ for all $n>k$ and the only non-zero elemets are $z_1,z_2,\dots,z_k$.\hfill$\square$

Let us remind some facts about Chevalley groups. We use definitions and de\-no\-ta\-tions from \cite{Car}.

Let $\Phi$ be an indecompasable root system of rang $l$ with the subsystem of simple roots $\Delta$, $|\Delta|=l$. The elementary Chevalley group $\Phi(F)$ of the type $\Phi$ over the field $F$ is a subgroup in the automorphism group of the simple Lie algebra $\mathcal{L}$  of the type $\Phi$, which is generated by the elementary root elements $x_{\alpha}(t)$, $\alpha\in\Phi$, $t\in F$. The dimension of the Lie algebra $\mathcal{L}$ is equal to $|\Phi|+|\Delta|$ and therefore group $\Phi(F)$ can be considered as a subgroup in the group of all $(|\Phi|+|\Delta|)\times(|\Phi|+|\Delta|)$ invertible matrices.

In the elementary Chevalley group, we consider the following important elements $n_{\alpha}(t)=x_{\alpha}(t)x_{-\alpha}(-t^{-1})x_{\alpha}(t)$,  $h_{\alpha}(t)=n_{\alpha}(t)n_{\alpha}(-1)$, $t\in F^*,~\alpha\in\Phi$.

For an arbitrary Chevalley group $G$ of the type $\Phi$ over the field $F$ we have the following short exact sequence of groups
$$1\to Z(G)\to G\to\Phi(F)\to1,$$
where $Z(G)$ is a center of the group $G$, and by the lemma \ref{pr3} we are mostly interested in the study of the $R_{\infty}$-property for elementary Chevalley groups.

Detailed information on the automorphism group of Chevalley groups can be found in \cite{Nas2, JHam}. Every Chevalley group has the following automorphisms
 \begin{enumerate}
 \item Inner automorphism $\varphi_g$, induced by the element $g\in G$
 $$\varphi_g:x\mapsto gxg^{-1}.$$
 \item Diagonal automorphism $\varphi_h$
 $$\varphi_h:x\mapsto hxh^{-1},$$
 where the element $h$ can be presented as a diagonal $(|\Phi|+|\Delta|)\times(|\Phi|+|\Delta|)$ matrix. If $F$ is an algebraically closed field then any diagonal automorphism is inner \cite[Lemma 4]{Nas2}.
 \item Field automorphism $\overline{\delta}$
 $$\overline{\delta}:x=(x_{ij})\mapsto (\delta(x_{ij})),$$
 where $\delta$ is an automorphism of the field $F$.
 \item Graph automorphism $\overline{\rho}$, which acts on the generators of the group $G$ by the rule
     $$\overline{\rho}:x_{\alpha}(t)\mapsto x_{\rho(\alpha)}(t),$$
     where $\rho$ is a symmetry of Dynkin diagram. An order of the graph automorphism is equal to 2 or to 3.
 \end{enumerate}

 Any field automorphism commutes with any graph automorphism. All the diagonal automorphisms form a normal subgroup in the group which is generated by diagonal, graph and field automorphisms.

 Theorem of Steinberg says that for any automorphism $\varphi$ of the elementary Chevalley group $G=\Phi(F)$ there exists an inner automorphism $\varphi_g$, a diagonal automorphism $\varphi_h$, a graph automorphism $\overline{\rho}$ and a field automorphism $\overline{\delta}$, such that $\varphi=\overline{\rho}\overline{\delta}\varphi_h\varphi_g$ \cite{JHam}.

\section{Chevalley groups}
In this paragraph we extend the following result from \cite[Theorem 1]{Nas2}.

\begin{ttt}\label{t1} Let $G$ be a Chevalley group of the type $\Phi$ over the field $F$ of zero characteristic and the transcendence degree of $F$ over $\mathbb{Q}$ is finite. Then
\begin{enumerate}
\item If $\Phi$ is a root system of the type $A_l(l\geq 7)$, $B_l(l\geq 4)$, $E_8$, $F_4$, $G_2$, then $G$ possesses the $R_{\infty}$-property.

\item If the equation $T^k=a$ can be solved in the field $F$ for any element $a$, then $G$ possesses the $R_{\infty}$-property in the case of the root systems $A_l(l=2,3,4,5,6)$, $B_l(l = 2, 3)$, $C_l(l\geq3)$, $D_l (l\geq 4)$, $E_6$, $E_7$, where $k$ is a positive integer from the table
    \begin{center}
\begin{tabular}[t]{|p{2em}|p{2em}|p{2em}|p{2em}|p{2em}|p{2em}|p{2em}|}
\hline
$\Phi$&~$A_l$  &~$B_l$ &$~C_l$&~$D_l$&~$E_6$&~$E_7$\\
\hline
k&$l+1$ &~~$2$&~~$2$&~~$2$&~~$6$&~~$2$ \\
\hline
\end{tabular}
\end{center}
\end{enumerate}
\end{ttt}

In particular, this theorem says that if $F$ is an algebraically closed field of zero characteristic, such that the transcendence degree of $F$ over $\mathbb{Q}$ is finite, then a Chevalley group of any normal type over the field $F$ possesses the $R_{\infty}$-property.

Here we exclude the condition of solvability of equations from the second item of the theorem \ref{t1}. We prove the following result.
\begin{ttt}\label{t2} Let $G$ be a Chevalley group of the type $\Phi$ over the field $F$ of zero characteristic. If the transcendence degree of $F$ over $\mathbb{Q}$ is finite, then $G$ possesses the $R_{\infty}$-property.
\end{ttt}
\textbf{Proof.}
Since $G/Z(G)\cong \Phi(F)$ then by the lemma \ref{pr3} it is sufficient to prove that the elementary Chevalley group  $\Phi(F)$ possesses the $R_{\infty}$-property. Consider that
 $G=\Phi(F)$.

 Let us consider an arbitrary automorphism $\varphi$ of the group $G$ and prove that the number of $\varphi$-conjugacy classes is infinite.
By the theorem of Steinberg there exists an inner automorphism $\varphi_g$, a diagonal automorphism $\varphi_h$, a graph automorphism $\overline{\rho}$ and a field automorphism $\overline{\delta}$, such that $\varphi=\overline{\rho}\overline{\delta}\varphi_h\varphi_g$. By the lemma \ref{pr1} the Reidemeister number $R(\varphi)$ is infinite if and only if the Reidemeister number $R(\varphi\varphi_{g^{-1}})$ is infinite, and we can consider that $\varphi=\overline{\rho}\overline{\delta}\varphi_{h}$.

Suppose that $R(\varphi)<\infty$ and consider the following elements of the group $G$
 $$g_i=h_{\alpha_1}(p_{i1})h_{\alpha_2}(p_{i2})\dots{h_{\alpha_l}}(p_{il}),~~~ i=1,2,\dots,$$
 where $p_{11}<p_{12}<\dots<p_{1l}<p_{21}<p_{22}<\dots$ are prime numbers. In the matrix representation the element $g_i$ has diagonal form
 $$g_i=diag(a_{i1},a_{i2},\dots,a_{i|\Phi|},\underbrace{1,\dots,1}_{|\Phi|}),$$
for certain rational numbers $a_{ij}$, such that $\nu(a_{ij})\neq\varnothing$ and $\nu(a_{ij})\cap\nu(a_{rs})=\varnothing$ for $i\neq r$ since $\nu(a_{ij})\subseteq\{p_{i1},\dots,p_{il}\}$ (see \cite{Nas2}).

Since $R(\varphi)<\infty$ then there exists an infinite subset of $\varphi$-conjugated elements in the set $g_1,g_2,\dots$. Without loosing of generality we can consider that all the elements $g_1,g_2,\dots$ belong to the $\varphi$-conjugacy class $[g_1]_{\varphi}$ of the element $g_1$.
It means that there exists an infinite set of matrices $Z_2,Z_3,\dots$ from $G$ such that
$$g_1=Z_ig_i\varphi(Z_i^{-1}),~~~i=2,3,\dots$$

Acting on this equalities by degrees of the automorphism $\varphi$ we have
\begin{eqnarray}\label{eq12}
\nonumber  g_1&=&Z_ig_i\varphi(Z_i^{-1}) ,\\
\nonumber  \varphi(g_1)&=&\varphi(Z_i)\varphi(g_i)\varphi^2(Z_i^{-1}),\\
\nonumber \varphi^2(g_1)&=&\varphi^2(Z_i)\varphi^2(g_i)\varphi^3(Z_i^{-1}),~~~i=2,3,\dots\\
\nonumber &\dots&\\
\nonumber \varphi^{5}(g_1)&=&\varphi^{5}(Z_i)\varphi^{5}(g_i)\varphi^{6}(Z_i^{-1}).
  \end{eqnarray}

If we multiply all of these equalities we conclude, that

\begin{equation}\label{main}
g_1\varphi(g_1)\dots\varphi^{5}(g_1)=Z_ig_i\varphi(g_i)\dots\varphi^{5}(g_i)\varphi^{6}(Z_i^{-1}) \end{equation}
Since the matrix $g_i$ has a diagonal form and the automorphism $\varphi_h$ acts as a conjugation by the diagonal matrix then $\varphi_h(g_i)=g_i$. Since the matrix $g_i$ has rational entries, then $\overline{\delta}(g_i)=g_i$ and therefore $\varphi(g_i)=\overline{\rho}(g_i)$.
If we denote $\widetilde{g}_i=g_i\varphi(g_i)\dots\varphi^{5}(g_i)=g_i\overline{\rho}(g_i)\dots\overline{\rho}^{5}(g_i)$
then  $$\widetilde{g}_i=diag(b_{i1},b_{i2},\dots,b_{i|\Phi|},\underbrace{1,\dots,1}_{l}),~~~i=1,2,\dots$$
since $\overline{\rho}$ permutes elements on the diagonal of the matrix $g_i$. Moreover, $\nu(b_{ij})\neq\varnothing$ and $\nu(b_{ij})\cap\nu(b_{rs})=\varnothing$ for $i\neq r$, since $\nu(b_{ij})\subseteq\nu(a_{i1})\cup\dots\cup\nu(a_{i|\Phi|})$.

 Since graph and field automorphisms commute and diagonal automorphisms form a normal subgroup in the group, which is generated by graph, field and diagonal automorphisms, then for a certain diagonal automorphism $\varphi_{\widetilde{h}}$ we have $\varphi^6=(\overline{\rho}\overline{\delta}\varphi_h)^6=\varphi_{\widetilde{h}}\overline{\delta}^6\overline{\rho}^6$. Since an order of the automorphism $\overline{\rho}$ is equal to $2$ or to $3$, then $\overline{\rho}^6=id$ and $\varphi^6=\varphi_{\widetilde{h}}\overline{\delta}^6$.
Then the equality (\ref{main}) can be rewritten
$$
\widetilde{g}_1=Z_i\widetilde{g}_i\varphi^{6}(Z_i^{-1})=Z_i\widetilde{g}_i\varphi_{\widetilde{h}}\overline{\delta}^6(Z_i^{-1})=Z_i\widetilde{g}_i{\widetilde{h}}\overline{\delta}^6(Z_i^{-1}){\widetilde{h}}^{-1},~~~ i=2,3,\dots
$$

If we multiply this equality by the element ${\widetilde{h}}$ on the right and denote $\hat{g}_i={\widetilde{g}}_i{\widetilde{h}}$ then we have
\begin{equation}\label{hren}
\hat{g}_1=Z_i\hat{g}_i\overline{\delta}^6(Z_i^{-1}),~~~ i=2,3,\dots
\end{equation}
From this equality we have
\begin{equation}\label{final}
\overline{\delta}^6(Z_i)=\hat{g}_1^{-1}Z_i\hat{g}_i,~~~i=2,3,\dots
\end{equation}
If we denote ${\widetilde{h}}=diag(c_{1},c_{2},\dots,c_{|\Phi|},\underbrace{1,\dots,1}_{l}),$ then $$\hat{g}_i=\widetilde{g}_i{\widetilde{h}}=diag(b_{i1}c_{1},b_{i2}c_{2},\dots,b_{i|\Phi|}c_{|\Phi|},\underbrace{1,\dots,1}_{l}),~~~~i=2,3,\dots$$
Let $\newcommand{\tempa}{\multicolumn{1}{c|}{Q_i}}
\newcommand{\tempb}{\multicolumn{1}{|c}{T_i}}
Z_i=\begin{pmatrix}
\tempa&R_i\\\cline{1-2}
S_i&\tempb
\end{pmatrix}$, where $Q_i=(q_{i,mn})$ is a  $|\Phi|\times|\Phi|$ matrix, $R_i=(r_{i,mn})$ is a $|\Phi|\times|\Delta|$ matrix, $S_i=(s_{i,mn})$ is a $|\Delta|\times|\Phi|$ matrix, $T_i=(t_{i,mn})$ is a  $|\Delta|\times|\Delta|$ matrix.

Then by the equality (\ref{final}) for all $m=1,\dots,|\Phi|,~n=1,\dots,|\Phi|$ we have
$$\delta^6(q_{i,mn})=(b_{1m}c_{m})^{-1}b_{in}c_{n}q_{i,mn}=d_{mn}b_{in}q_{i,mn},~~~i=2,3,\dots,$$
where $d_{mn}=(b_{1m}c_{m})^{-1}c_{n}$. Since $\nu(b_{in})\neq\varnothing$ and $\nu(b_{in})\cap\nu(b_{jn})=\varnothing$ for $i\neq j$, then we can apply the lemma \ref{lem2} to the set $q_{2,mn}, q_{3,mn},\dots$ Therefore by the lemma \ref{lem2} there exists a positive integer $N_{mn}$, such that $q_{i,mn}=0$ for every $i>N_{mn}$.

If we denote by $N$ the value
$$N=\max_{
n,m=1,\dots,|\Phi|} N_{mn},$$
then for every $i>n$ we have $Q_i=O_{|\Phi|\times|\Phi|}$. Using the same arguments to the matrices $S_2,S_3,\dots$ we conclude that for sufficiently large indexes $i$ all the matrices $S_i$ are the matrices with zero entries only, and therefore the matrix $Z_i$ has the form $$Z_i=\begin{pmatrix}O_{|\Phi|\times|\Phi|}&R_i\\
O_{|\Delta|\times|\Phi|}&T_i
\end{pmatrix}.$$
The determinant of this matrix is equal to zero and it means that $Z_i$ can not belong to $G$. This contradiction proves the theorem.
\begin{flushright}
$\square$
\end{flushright}

\section{Linear algebraic groups}

If $G$ is a linear algebraic group over an algebraically closed field, then it has a unique maximal solvable normal subgroup $R(G)$, called the radical of $G$. A connected linear algebraic group $G$ is called reductive if its radical is a torus, or, equivalently, if it can be decomposed $G=G^{\prime}T^{\prime}$ with $G^{\prime}$ a semisimple group and $T^{\prime}$ a central torus \cite[\S6.5]{S}.

The quotient group $G/R(G)$ has a trivial radical, i.e. is a semisimple group \cite[\S 19.5]{J}.

\begin{ttt}\label{tema} Let $F$ be such an algebraically closed field of zero characteristic that the transcendence degree of $F$ over $\mathbb{Q}$ is finite. If the reductive linear algebraic group $G$ over the field $F$ has a nontrivial quotient group $G/R(G)$, then $G$ possesses the $R_{\infty}$-property.
\end{ttt}
\textbf{Proof.} For the group $G$ we have the following short exact sequence of groups
$$1\to R(G)\to G\to G/R(G)\to1,$$
Since $G$ is reductive, then the radical $R(G)$ is a central torus and therefore is  a characteristic subgroup of $G$. Hence by the lemma \ref{pr3} it is sufficient to prove that the semisimple group $G/R(G)$ possesses the $R_{\infty}$-property and we can consider that $G$ is a semisimple linear algebraic group. Since $G$ is a semisimple linear algebraic group, then it is a product, with some amalgamation of (finite) centers, of its simple subgroups $H_1, H_2,\dots, H_k$ \cite[\S14.2]{J}
$$G=H_1\circ\dots \circ H_k.$$
Every simple linear algebraic group $H_i$ is a Chevalley group of (normal) type $\Phi_i$ over the field $F$. Factoring the group $G$ by its center we have the following short exact sequence of groups
 $$1\to Z(H_1\circ\dots \circ H_k)\to H_1\circ\dots \circ H_k\to \Phi_1(F)\times\dots\times\Phi_k(F)\to1,$$
 where $\Phi_i(F)$ is an elementary Chevalley group of the type $\Phi_i$ over the field $F$. Hence, by the lemma \ref{pr3} we can consider that $G=\Phi_1(F)\times\dots\times\Phi_k(F)$ and prove that this group possesses the $R_{\infty}$-property. Permute the groups $\Phi_1(F),\dots, \Phi_k(F)$ so that all the groups with the same root system form blocks
 $$G=\underbrace{\Phi_1(F)\times\dots\Phi_1(F)}_{k_1}\times\underbrace{\Phi_2(F)\times\dots\Phi_2(F)}_{k_2}\times\dots\times\underbrace{\Phi_r(F)\times\dots\Phi_r(F)}_{k_r},$$
 where $k_1+k_2+\dots+k_r=k$. Denote $G_i=\underbrace{\Phi_i(F)\times\dots\Phi_i(F)}_{k_i}$.

 Every group $G_i$ is a characteristic subgroup of  $G=G_1\times\dots\times G_r$. It is obvious that if some group is a direct product of its characteristic subgroups and at least one of this subgroups possesses the $R_{\infty}$-property then the group itself possesses the $R_{\infty}$-property. Therefore it is sufficient to prove that the group $$G=\Phi(F)\times\dots\times\Phi(F)=\Phi(F)^k$$ possesses the $R_{\infty}$-property.

Every element $g\in G=\Phi(F)^k$ can be presented as a direct sum of $k$ matrices $g_1,\dots,g_k$ of the size $(|\Phi|+|\Delta|)\times(|\Phi|+|\Delta|)$ each of which belongs to $\Phi(F)$.

An automorphism group of $G$ has the form
 \begin{equation}\label{Aut}{\rm Aut}(G)=({\rm Aut}(\Phi(F)))^k\leftthreetimes S_k,\end{equation}
where $S_k$ is a permutation group on $k$ symbols.

 To prove that the group $G=\Phi(F)^k$ possesses the $R_{\infty}$-property consider an arbitrary automorphism  $\varphi$ of the group $G$ and prove that $R(\varphi)=\infty$.
By the equality (\ref{Aut}) the automorphism $\varphi$ can be written in the following form
$$\varphi=\left(\varphi_1,\dots,\varphi_k,\sigma\right),$$
where $\varphi_1,\dots,\varphi_k\in {\rm Aut(\Phi(F))}$, $\sigma\in S_k$, and $\varphi$ acts on the group $G$ by the rule
\begin{equation}\label{prav}
\varphi: x_1\oplus x_2\oplus\dots\oplus x_k \mapsto \varphi_{1^{\sigma}}\left(x_{1^{\sigma}}\right)\oplus\varphi_{2^{\sigma}}\left(x_{2^{\sigma}}\right)\oplus\dots\oplus\varphi_{k^{\sigma}}\left(x_{k^{\sigma}}\right),
\end{equation}
where $i^{\sigma}$ denotes an image of $i$ by the permutation $\sigma$.

Every automorphism $\varphi_i\in {\rm Aut}(\Phi(F))$ can be presented as a product of the inner automorphism $\varphi_{g_i}$, the diagonal automorphism $\varphi_{h_i}$, the graph automorphism $\overline{\rho}_i$ and the field automorphism $\overline{\delta}_i$. Since $F$ is an algebraically closed field, then every diagonal automorphism $\varphi_{h_i}$ is inner \cite[Lemma 4]{Nas2}, hence for every $i$ we can consider that  $\varphi_i=\varphi_{x_i}\overline{\rho}_i\overline{\delta}_i$. Then the automorphism  $\varphi$ can be presented as a product of two automorphism
$$\varphi=(\varphi_{x_{1^{\sigma}}},\varphi_{x_{2^{\sigma}}},\dots,\varphi_{x_{k^{\sigma}}},id)(\overline{\rho}_1\overline{\delta}_1,\overline{\rho}_2\overline{\delta}_2,\dots,\overline{\rho}_k\overline{\delta}_k,\sigma),$$
where $(\varphi_{x_{1^{\sigma}}},\varphi_{x_{2^{\sigma}}},\dots,\varphi_{x_{k^{\sigma}}},id)$ is an inner automorphism. By the lemma \ref{pr1} we can consider that $\varphi_i=\overline{\rho}_i\overline{\delta}_i$ and
$$
\varphi=(\overline{\rho}_1\overline{\delta}_1,\overline{\rho}_2\overline{\delta}_2,\dots,\overline{\rho}_k\overline{\delta}_k,\sigma).
$$

Using the induction on $r$ prove that
\begin{equation}\label{step}\varphi^r: g_1\oplus\dots\oplus g_k\mapsto\psi_{1}(x_{1^{\sigma^r}})\oplus\dots\oplus\psi_{k}(x_{k^{\sigma^r}}),
\end{equation}
where $\psi_i=\varphi_{i^{\sigma}}\varphi_{i^{\sigma^2}}\dots\varphi_{i^{\sigma^r}}$.

The basis of the induction ($r=1$) is obvious (the equality (\ref{prav})). If we suppose that the equality (\ref{step}) holds for some  $r$, then
 \begin{eqnarray}
\nonumber \varphi^{r+1}(g_1\oplus\dots\oplus g_k)&=&\varphi(\varphi^r(g_1\oplus\dots\oplus g_k))\\
\nonumber &=&\varphi(\psi_{1}(x_{1^{\sigma^r}})\oplus\dots\oplus\psi_{k}(x_{k^{\sigma^r}}))\\
\nonumber &=&\varphi_{1^{\sigma}}\psi_{1^{\sigma}}(x_{1^{\sigma^{r+1}}})\oplus\dots\oplus\varphi_{k^{\sigma}}\psi_{k^{\sigma}}(x_{k^{\sigma^{r+1}}}).\end{eqnarray}
Note the equality $\varphi_{i^{\sigma}}\psi_{i^{\sigma}}=\varphi_{i^{\sigma}}\varphi_{i^{\sigma^2}}\dots\varphi_{i^{\sigma^{r+1}}}$, what we wanted to prove.

Consider the set of elements $g_1,g_2,\dots$ of the group $\Phi(F)$ from the theorem \ref{t2}
 $$g_i=h_{\alpha_1}(p_{i1})h_{\alpha_2}(p_{i2})\dots{h_{\alpha_l}}(p_{il}),~~~ i=1,2,\dots,$$
 where $p_{11}<p_{12}<\dots<p_{1l}<p_{21}<p_{22}<\dots$ are prime integers. This elements are presented by the diagonal matrices
 $$g_i=diag(a_{i1},a_{i2},\dots,a_{i|\Phi|},\underbrace{1,\dots,1}_{l}),~~~i=1,2,\dots$$
where $a_{ij}$ are rational numbers, such that $\nu(a_{ij})\neq\varnothing$ and $\nu(a_{ij})\cap\nu(a_{rs})=\varnothing$ for $i\neq r$.

As we already shown in the theorem \ref{t2} for every automorphism $\varphi_j=\overline{\rho}_j\overline{\delta}_j$ we have $\varphi_j(g_i)=\overline{\rho}_j(g_i)$.

Let us consider the set of elements $\widetilde{g}_1, \widetilde{g}_2,\dots$ of the group $G=\Phi(F)^k$, where $\widetilde{g}_i=g_i\oplus\dots\oplus g_i$. Then by the arguments above $\varphi(\widetilde{g}_i)=\overline{\rho}_{1^{\sigma}}(g_i)\oplus\dots\oplus \overline{\rho}_{k^{\sigma}}(g_i)$.

Suppose that $R(\varphi)<\infty$. Then there is an infinite subset of $\varphi$-conjugated elements in the set $\widetilde{g}_1,\widetilde{g}_2,\dots$. Without loosing of generality we can consider that all the matrices $\widetilde{g}_1,\widetilde{g}_2,\dots$ belong to the $\varphi$-conjugacy class $[\widetilde{g}_1]_{\varphi}$ of the element $\widetilde{g}_1$. Then for certain matrices $Z_2, Z_3,\dots$ we have
$$\widetilde{g}_1=Z_i\widetilde{g}_i\varphi(Z_i^{-1}),~~~i=2,3,\dots$$
Denote by $s$ an order of the permutation $\sigma$ and act on this equality by the degrees of the automorphism $\varphi$
\begin{eqnarray}
\nonumber \widetilde{g}_1&=&Z_i\widetilde{g}_i\varphi(Z_i^{-1})\\
\varphi(\nonumber \widetilde{g}_1)&=&\varphi(Z_i)\varphi(\widetilde{g}_i)\varphi^2(Z_i^{-1})\\
\nonumber\vdots&\vdots&\vdots\\
\nonumber \varphi^{6s-2}(\widetilde{g}_1)&=&\varphi^{6s-2}(Z_i)\varphi^{6s-2}(\widetilde{g}_i)\varphi^{6s-1}(Z_i^{-1})\\
\nonumber \varphi^{6s-1}(\widetilde{g}_1)&=&\varphi^{6s-1}(Z_i)\varphi^{6s-1}(\widetilde{g}_i)\varphi^{6s}(Z_i^{-1})
\end{eqnarray}

If we multiply all of these equalities, we obtain the  equality

\begin{equation}\label{last}
\widetilde{g}_1\varphi(\widetilde{g}_1)\varphi^2(\widetilde{g}_1)\dots\varphi^{6s-1}(\widetilde{g}_1)=Z_i\widetilde{g}_i\varphi(\widetilde{g}_i)\varphi^2(\widetilde{g}_i)\dots\varphi^{6s-1}(\widetilde{g}_i)\varphi^{6s}(Z_i^{-1}).
\end{equation}
The element $\widetilde{g}_i\varphi(\widetilde{g}_i)\varphi^2(\widetilde{g}_i)\dots\varphi^{6s-1}(\widetilde{g}_i)$ can be rewritten in details
\begin{eqnarray} \nonumber &&\widetilde{g}_i\varphi(\widetilde{g}_i)\varphi^2(\widetilde{g}_i)\dots\varphi^{6s-1}(\widetilde{g}_i)=\\
\nonumber&&=(g_i\oplus\dots\oplus g_i)(\overline{\rho}_{1^{\sigma}}(g_i)\oplus\dots\oplus \overline{\rho}_{k^{\sigma}}(g_i))\dots(\overline{\rho}_{1^{\sigma^{6s-1}}}(g_i)\oplus\dots\oplus \overline{\rho}_{k^{\sigma^{6s-1}}}(g_i))=\\ \nonumber &&=g_i\overline{\rho}_{1^{\sigma}}(g_i)\dots\overline{\rho}_{1^{\sigma^{6s-1}}}(g_i)\oplus\dots\oplus g_i\overline{\rho}_{k^{\sigma}}(g_i)\dots\overline{\rho}_{k^{\sigma^{6s-1}}}(g_i))=\widehat{g}_{i1}\oplus\dots\oplus \widehat{g}_{ik},
\end{eqnarray}
where $\widehat{g}_{ij}=g_i\overline{\rho}_{j^{\sigma}}(g_i)\dots\overline{\rho}_{j^{\sigma^{6s-1}}}(g_i)$.

Since every graph automorphism  $\overline{\rho}_j$ permutes elements on the diagonal of the matrix $g_i$ then for every $j=1,\dots,k$, $i=1,2,\dots$ we have
\begin{equation}
\widehat{g}_{ij}=diag(b_{ij1},b_{ij2},\dots,b_{ij|\Phi|},\underbrace{1,\dots,1}_{l}),
\end{equation}
where, $\nu(b_{ijr})\neq\varnothing$ and $\nu(b_{ijr})\cap\nu(b_{uvw})=\varnothing$ for $i\neq u$ since $\nu(b_{ijr})\subseteq\{p_{i1},\dots,p_{il}\}$.

By the equality (\ref{step}) we have $\varphi^s=(\psi_1,\psi_2,\dots,\psi_k,id)$, where $\psi_i=\varphi_{i^{\sigma}}\varphi_{i^{\sigma^2}}\dots\varphi_{i^{\sigma^r}}$. Since all the automorphisms $\varphi_1,\varphi_2,\dots,\varphi_k$ are the products of graph and field automorphism ($\varphi_i=\overline{\rho}_i\overline{\delta}_i$) and graph and field automorphisms commute then every automorphism $\psi_i$ is a product of graph and field automorphisms $\psi_i=\overline{\xi}_i\overline{\theta}_i$ for certain $\overline{\xi}_i$, $\overline{\theta}_i$. Therefore
$$\varphi^{6s}=(\varphi^s)^6=(\overline{\xi}_1\overline{\theta}_1,\dots \overline{\xi}_k\overline{\theta}_k,id)^6=(\overline{\xi}_1^6\overline{\theta}_1^6,\dots \overline{\xi}_k^6\overline{\theta}_k^6,id)=(\overline{\theta}_1^6,\dots \overline{\theta}_k^6,id)$$
Using this fact, denoting the matrix $Z_i=Z_{i1}\oplus\dots \oplus Z_{ik}$ projecting the equality (\ref{last}) to the first group $\Phi(F)$ we obtain the equality
$$\widehat{g}_{11}=Z_{i1}\widehat{g}_{i1}\overline{\theta}_1^6(Z_{i1}),~~~i=2,3,\dots$$
  This equality is the same as the equality (\ref{hren}) from the theorem \ref{t2}. Using the same arguments as in the theorem \ref{t2} we conclude that for sufficiently large coefficient $N$ the matrix  $Z_{iN}$ is degenerated and therefore the matrix $Z_{N}$ is degenerated but it contradicts to the fact that $Z_N$ belongs to $G$.\hfill$\square$

We use the fact that the group $G$ is a reductive linear algebraic group in order to say that the radical $R(G)$ is a characteristic subgroup of $G$. Even the theorem \ref{tema} holds for every connected linear algebraic group such that the radical $R(G)$ is a characteristic. For example, if any automorphism of the group $G$ is a morphism of the group $G$ (as of an affine manifold) then the radical $R(G)$ is characteristic \cite[Theorem 7.1(c)]{S} and the theorem \ref{tema} holds for such groups.

\section{Finite Reidemeister number in linear groups}

Following \cite{Romankov}, we define the {\it Reidemeister spectrum
of} $G$ as
$$Spec(G)=\{ R(\varphi) \: | \: \varphi \in \Aut(G)\}.$$
In particular, $G$ possesses the $R_\infty$-property if and only if $Spec(G)=\{\infty\}$.

 It is easy to see that
  $Spec({\mathbb Z}) = \{2\} \cup \{ \infty \}$,
and,  for $n \geq 2$, the spectrum of ${\mathbb
Z}^n$ is full, i.e. $ Spec({\mathbb
Z}^n) = {\mathbb N} \cup \{\infty \}$. For free nilpotent groups
we have the following: $Spec(N_{22}) = 2{\mathbb N} \cup \{\infty
\}$ ($N_{22}$ is the discrete Heisenberg group)
\cite{Indukaev,FelIndTro,Romankov},  $ Spec(N_{23}) = \{ 2k^2 \: | \: k
\in {\mathbb N}\} \cup \{\infty \}$ \cite{Romankov} and $ Spec
(N_{32}) = \{2k - 1\: | \: k \in {\mathbb N}\} \cup \{ 4k \: | \: k \in
{\mathbb N}\} \cup \{ \infty \}$ \cite{Romankov}.

Recently, in \cite{DeGon} it was proven that the group
$N_{rc}~(r > 1)$ admits an automorphism with finite Reidemeister number if and only if $c < 2r$.

 In \cite{GoWon},  examples of polycyclic non-virtually nilpotent groups  which admit automorphisms with
 finite Reidemeister numbers have been described. In this examples $G$ is a semidirect
product of $\Z^2$ and $\Z$ by Anosov automorphism defined by the
matrix $\scriptsize\Mat 2111$. The group $G$ is solvable and of the exponential
growth. The automorphism $\varphi$ with finite Reidemeister number is defined by $\scriptsize\Mat
{~~0}1{-1}0$ on $\Z^2$ and as $-id$ on $\Z$.

 Metabelian (therefore, solvable) finitely generated, non-polycyclic groups have quite interesting
 Reidemeister spectrum \cite{FelGon2011Q}:
  for example, if the homomorphism $\theta:\mathbb Z\to{\rm Aut}(\mathbb Z[1/p]^2)$ is such that $\theta(1)=\scriptsize\Mat r00s,$ then we have the following cases:
\begin{itemize}
\item[a)] If  $r=s=\pm 1$ then $Spec(\mathbb Z[1/p]^2\rtimes_{\theta}
\mathbb Z)=\{2n\: | \:  n\in \mathbb N, \ (n,p)=1    \}\cup \{\infty\}$,
where $(n,p)$   denote the greatest common divisor of $n$ and $p$.

\item[b)] If $r=-s=\pm 1$ then  $Spec( \mathbb
Z[1/p]^2\rtimes_{\theta}\mathbb Z)=\{ 2p^l(p^k\pm 1), 4p^l \: | \: l,
k>0\} \cup\{  \infty\}$.

\item[c)] If  $rs=1$ and $|r|\ne 1$ then   $Spec(\mathbb
Z[1/p]^2\rtimes_{\theta}\mathbb Z)=\{ 2(p^l \pm 1), 4 \: | \:l>0 \}
\cup\{ \infty \}$.

\item[d)] If either $r$ or $s$
does not equal to $\pm1$,
and $ rs \ne 1$ then
  $Spec(\mathbb Z[1/p]^2\rtimes_{\theta}\mathbb
Z)=\{\infty\}$.
\end{itemize}
In the paper \cite{Jab} Jabara proved that if residually finite group $G$ admits an automorphism of prime order $p$ with finite Reidemeister number, then $G$ is virtually nilpotent group of class bounded by a function of $p$.

From another side, we have described in an Introduction a lot of classes of non- solvable, finitely generated, residually finite groups which have the $R_{\infty}$-property.
All together  was a motivation for  the following conjecture \\
~\\
\textbf{Conjecture R (A. Fel'shtyn, E. Troitsky \cite[Conjecture R]{FTMPI})} Every infinite, residually finite, finitely generated group either possesses the $R_{\infty}$-property or is a virtually solvable group.\\

 Here we study this question for infinite linear groups.
\begin{prp}\label{prpr}
Let $G$ be a reductive linear algebraic over the field $F$ of zero characteristic and finite transcendence degree over $\mathbb{Q}$. If $G$ possesses an automorphism $\varphi$ with finite Reidemeister number then $G$ is a torus.
\end{prp}
\textbf{Proof} Since $G$ possesses an automorphism $\varphi$ with finite Reidemeister number,
then by the theorem \ref{tema}, it has trivial quotient group $G/R(G)$, therefore $G=R(G)$ and
hence $G$ is a central torus (therefore, is solvable).\hfill$\square$

%

\section{Groups with property $S_\infty$}\label{sec:isogred}


Suppose that $\Psi \in \Out (G)= \Aut (G)/ \Inn (G)$.
Let $\cS(\Psi)$ be the set of isogredience classes of $\Psi$.
Then  $\cS(\Id)$ can be identified with the set of conjugacy classes
of  $G/Z(G)$ (see \cite{FTMPI}).

 The
definition of the similarity (isogredience) from Introduction  goes back to Jacob  Nielsen.
 He observed (see \cite{Jiang})
 that conjugate lifting of homeomorphism of
 surface  have similar dynamical properties.
 This led Nielsen to the definition of the
 isogredience of liftings in this case.
 Later Reidemeister  and
 Wecken succeeded in generalizing the theory to
 continuous maps of compact polyhedra
 (see \cite{Jiang}).

 The set of isogredience classes of
 automorphisms representing a
 given outer automorphism and the
 notion of index $\Ind(\Psi)$
 defined via the set of isogredience
 classes are strongly related to important
 structural properties of $\Psi$
 (see \cite{GabJaegLevittLustig98Duke}).

One of the main results of \cite{ll} is that for any
non-elementary hyperbolic group and any $\Psi$ the set
$ \cS(\Psi)$ is infinite, i.~e.
$S(\Psi)=\infty$.
Thus, this result  says:
any non-elementary hyperbolic group is an $S_\infty$-group.
On the other hand, finite and finitely generated Abelian groups are
evidently non $S_\infty$-groups.

Two representatives of $\Psi$ have the forms $\varphi_sa$,
$\varphi_qa$ for some $s,q \in G$ and fixed $a \in \Psi$. They are isogredient if and only
if
$$
\varphi_qa =\varphi_g\varphi_s  a\varphi_g^{-1}
=\varphi_g  \varphi_s\varphi_{a(g^{-1})}a,
$$
$$
\varphi_q=\varphi_{gsa(g^{-1})},\qquad q=gsa(g^{-1})c,\quad c\in
Z(G)
$$
(see \cite[p. 512]{ll}).
So, the following statement is proved.
\begin{lem}\cite[Lemma 3.3]{FTMPI}\label{lem:SandR} Let $\varphi\in\Psi$ be an automorphism of the group $G$ and $\overline{\varphi}$ be an automorphism of the group $G/Z(G)$ which is induced by $\varphi$. Then the number $S(\Psi)$ is equal to the number of $\overline{\varphi}$-conjugacy classes in the group $G/Z(G)$.
\end{lem}

Since $Z(G)$ is a characteristic subgroup, we obtain
 the following statement

\begin{ttt}\cite[Theorem 3.4]{FTMPI}\label{teo:SandRfinZ}
Suppose,
$|Z(G)|<\infty$.
Then $G$ is an $R_\infty$-group
if and only if
$G$ is an $S_\infty$-group.
\end{ttt}

 A more advanced example
of a non $S_\infty$-group is the
Osin's group \cite{Osin}. This is a non-residually
finite exponential growth group with two conjugacy
classes. Since it is simple, it is not $S_\infty$
group (see \cite{FTMPI}).

\begin{ttt}\label{sim} Let $F$ be such an algebraically closed field of zero characteristic that the transcendence degree of $F$ over $\mathbb{Q}$ is finite. If the reductive  linear algebraic group $G$ over the field $F$ has a nontrivial quotient group $G/R(G)$, then $G$ possesses the $S_{\infty}$-property.
\end{ttt}
\textbf{Proof.} Since $R(G)$ is a characteristic subgroup of $G$ then by the lemma \ref{an3} it is sufficient to prove the theorem for semisimple group $G/R(G)$. The result follows immediately from the theorem \ref{tema} and the theorem \ref{teo:SandRfinZ} and from the fact that semisimple linear algebraic group has  finite center.
\hfill$\square$

\begin{prp}\label{prpr2}
Let $G$ be a reductive linear algebraic group over the field $F$ of zero characteristic and finite transcendence degree over $\mathbb{Q}$. If $G$ possesses an outer automorphism $\Psi$ with finite number $S(\Psi)$ then $G$ is a torus.
\end{prp}
\textbf{Proof}  Since $G$  possesses an outer automorphism  $\Psi$ with finite  number $S(\Psi)$, then by the theorem \ref{sim}, it has a trivial quotient  $G/R(G)$, therefore $G=R(G)$ and is a central torus.
\hfill$\square$


The conjecture  of Fel'shtyn and Troitsky from the section 5 can be rewritten in terms of $S_{\infty}$-property by the following way.\\
~\\
\textbf{Conjecture S} Every infinite, residually finite, finitely generated group either possesses the $S_{\infty}$-property or is a virtually solvable group.\\

Really, if $S(\varphi {\rm Inn}(G))<\infty$ for some automorphism $\varphi\in {\rm Aut}(G)$ then by the lemma \ref{lem:SandR} we have $R(\overline{\varphi})<\infty$, where $\overline{\varphi}$ is an automorphism of the group $G/Z(G)$ induced by $\varphi$. Since $G$ is residually finite finitely generated group, then $G/Z(G)$ is also finitely generated and residually finite and by the conjecture R is a virtually solvable group.

It means that there exists a solvable subgroup $\overline{H}\leq G/Z(G)$ of finite index. Let $n$ be a derived length of $\overline{H}$, i.~e. $\overline{H}^{(n)}=1$. Let $H$ be a preimage of $\overline{H}$ under the canonical homomorphism $G\to G/Z(G)$. Then $H^{(n)}\leq Z(G)$ and $H^{(n+1)}=1$, therefore $H$ is a solvable group. Since $G/H\simeq (G/Z(G))/(H/Z(G))=(G/Z(G))/\overline{H}$, then the index of $H$ in $G$ is equal to the index of $\overline{H}$ in $G/Z(G)$, i.~e. is finite, therefore $G$ is solvable by finite.

We have proven in all that \cite[Conjecture S]{FTMPI} can be formulated without
the restriction  that  the group under consideration  has finite center.

~\\
 {\scshape Instytut Matematyki, Uniwersytet Szczecinski, ul. Wielkopolska 15, 70-451 Szczecin,
Poland}\\
\emph{E-mail address:} fels@wmf.univ.szczecin.pl

~\\
 {\scshape Sobolev Institute of mathematics, ak. Koptyug avenue 4, 630090, Novosibirsk, Russia}\\
\emph{E-mail address:} ntr@math.nsc.ru, timur.nasybullov@mail.ru

\end{document}